\documentclass[proc]{edpsmath}

\usepackage{hyperref}

\usepackage{amsfonts,amssymb,stmaryrd,amsmath,amssymb}
\usepackage{graphicx,color}
\usepackage{multimedia}
\usepackage{float}
\usepackage{placeins}

\usepackage{moreverb}

\newtheorem{lemma}{Lemma}

\newcommand \bn {\boldsymbol{\mathrm{n}}}

\newcommand \bh {\boldsymbol{\mathrm{h}}}
\newcommand \bomega {\boldsymbol{\mathrm{\omega}}}

\newcommand \bu {\boldsymbol{\mathrm{u}}}
\newcommand \bv {\boldsymbol{\mathrm{v}}}
\newcommand \bg {\boldsymbol{\mathrm{g}}}
\newcommand \bff {\boldsymbol{\mathrm{f}}}
\newcommand \blambda {\boldsymbol{\mathrm{\lambda}}}
\newcommand \bmu {\boldsymbol{\mathrm{\mu}}}

\newcommand \p {\partial}

\newcommand \x {\mathrm{x}}
\newcommand \y {\mathrm{y}}

\newcommand \R {\mathbb{R}}

\renewcommand \L {\mathrm{L}}

\renewcommand \H {\mathrm{H}}

\renewcommand \d {\mathrm{d}}

\renewcommand \div {\mathrm{div}}

\begingroup\makeatletter\ifx\SetFigFont\undefined%
\gdef\SetFigFont#1#2#3#4#5{%
  \reset@font\fontsize{#1}{#2pt}%
  \fontfamily{#3}\fontseries{#4}\fontshape{#5}%
  \selectfont}%
\fi\endgroup%

\begin{document}

%%-----------------------------
%%      the top matter
%%-----------------------------
\title{A fictitious domain approach for Fluid-Structure Interactions based on the eXtended Finite Element Method.}\thanks{This work is partially supported by the ANR-project CISIFS: 09-BLAN-0213-03, and the foundation STAE in the context of the RTRA platform SMARTWING. It is in particular based on the use of the Getfem++ library \cite{Getfem}, and thus on collaborative efforts with Yves Renard.}% At most 5 thanks
\author{S\'ebastien COURT}\address{Universit\'e Blaise Pascal, Laboratoire de Math\'ematiques, Campus des C\'ezeaux, B.P. 80026, 63171 Aubi\`ere Cedex France;\\ E-mail: \href{mailto:sebastien.court@math.univ-bpclermont.fr}{\texttt{sebastien.court@math.univ-bpclermont.fr}}}
\author{Michel FOURNI\'E}\address{Universit\'e de Toulouse, Institut de Math\'ematiques de Toulouse, 118 route de Narbonne, 31062 TOULOUSE Cedex;\\ E-mail: \href{mailto:michel.fournie@math.univ-toulouse.fr}{\texttt{michel.fournie@math.univ-toulouse.fr}}}
\author{Alexei LOZINSKI}\address{Universit\'e de Franche-Comt\'e, Laboratoire de Math\'ematiques de Besan\c con, 16 route de Gray, 25030 Besançon Cedex;\\ E-mail: \href{mailto:alexei.lozinski@univ-fcomte.fr}{\texttt{alexei.lozinski@univ-fcomte.fr}}}
%
%\dedicated{\it Dedicated to Maurice Dupont} %if necessary
%
\begin{abstract}
In this work we develop a {\it fictitious domain} method for the Stokes problem which allows computations in domains whose boundaries do not depend on the mesh. The method is based on the ideas of Xfem and has been first introduced for the Poisson problem. The fluid part is treated by a mixed finite element method, and a Dirichlet condition is imposed by a Lagrange multiplier on an immersed structure localized by a {\it level-set} function. A stabilization technique is carried out in order to get the convergence for this multiplier. The latter represents the forces that the fluid applies on the structure. The aim is to perform fluid-structure simulations for which these forces have a central role. We illustrate the capacities of the method by extending it to the incompressible Navier-Stokes equations coupled with a moving rigid solid.
\end{abstract}
\begin{resume}
Dans ce travail nous d\'eveloppons une m\'ethode de type {\it domaines fictifs} pour le probl\`eme de Stokes permettant des calculs dans des domaines ind\'ependants du maillage. La m\'ethode est bas\'ee sur les \textcolor{black}{id\'ees} de Xfem et a d'abord \'et\'e introduite pour le probl\`eme de Poisson. La partie fluide est trait\'ee avec une m\'ethode \'el\'ements finis mixte, et une condition de Dirichlet est impos\'ee par un multiplicateur de Lagrange sur une structure immerg\'ee localis\'ee par une fonction {\it level-set}. \textcolor{black}{Une} technique de stabilisation est mise en oeuvre pour obtenir la convergence de ce multiplicateur. Ce dernier repr\'esente les forces que le fluide exerce sur la structure. L'objectif est de r\'ealiser des simulations fluide-structure pour \textcolor{black}{lesquelles} ces forces ont un r\^ole central. Nous illustrons les capacit\'es de la m\'ethode en l'\'etendant aux \'equations de Navier-Stokes \textcolor{black}{incompressibles} coupl\'ees avec un solide rigide en mouvement.
\end{resume}
\maketitle
%%-----------------------------
%%      your text
%%-----------------------------
\section*{Introduction}
The development of efficient numerical methods for the simulation of fluid-structure interaction problems is a great challenge. On one hand some methods are based on meshes which fit to the geometry of the computational domain, like in \cite{LT, SMSTT0, SST} for instance. In that case re-meshing is necessary when the geometry is led to evolve, which can be greedy in resources, especially for complex problems. On the other hand numerical methods consider meshes with a fictitious domain approach where the mesh is cut by boundaries. Let us cite the Immersed Boundary Method (see \cite{Peskinacta}) for which force terms are added on a boundary representing the fluid-structure interface. Let us also mention the Lagrange multiplier method, carried out for the displacement of rigid solids in incompressible fluids (see \cite{Glowinski}). In this approach the velocity of the fluid is extended inside the solid by taking the values of the velocity of the latter. More recently the eXtended Finite Element Method (Xfem in abbreviated form) developed by Mo\"es, Dolbow and Belytschko in \cite{Moes} for cracked domains problems (see \cite{MoesB, SukumarC, Gerstenberger2008, Choi2010} for instance) has been applied to fluid-structure interaction problems. It is based on a representation by {\it level-set} of the boundary (originally the crack) and on the enrichment of a finite element space by singular functions near the boundary.\\
In the framework of fluid-structure interactions, the difficulty that present the techniques aforementioned lies in the choice of the Lagrange multiplier associated with the boundary conditions taken into account at the interface, because of the fact that the latter cuts the mesh. Moreover, our concern requires that this multiplier space has to satisfy an inf-sup condition which leads to a good approximation of the multiplier. Indeed, the multiplier associated with the Dirichlet condition (which is imposed by the equality of velocities) is the normal trace of the Cauchy stress tensor $\sigma(\bu, p)\bn$ at the interface, and represents nothing else the force that the fluid exerts on the structure. Thus getting the convergence of this quantity is a crucial point. The method we develop has been first introduced in \cite{HaslR} for the Poisson problem, and then adapted to the Stokes problem in \cite{Court}, which is the keystone of problems involving viscous incompressible fluids, notably the full incompressible Navier-Stokes equations. It is based on the ideas of Xfem, but it combines a fictitious domain approach with a stabilization technique for the convergence of $\sigma(\bu,p)\bn$. Note that the alternative methods based on the work of Nitsche \cite{Nitsche} (see \cite{Burman1, Burman3, Massing} for instance) do not introduce the Lagrange multiplier and so do not necessarily provide a good approximation of this force.\\ 
The paper is organized as follows: In section \ref{sec1} we present the Stokes problem we consider and the extended variational formulation associated with the search of a weak solution to this system. Then in section \ref{sec3} we describe the fictitious domain approach we develop and we give a theoretical result dealing with the convergence for the complete method. The robustness of the method with respect to the geometrical configurations is highlighted in section \ref{sec3}. Section \ref{sec4} is devoted to more complex numerical simulations involving the Navier-Stokes equations coupled with a moving rigid solid, before conclusion.

\section{Framework} \label{sec1}
Given a bounded domain $\mathcal{O} \subset \R^2$, we consider a solid immersed in a viscous incompressible fluid and occupying a domain denoted by $\mathcal{S}$ (with boundary denoted by $\Gamma$). The surrounding fluid occupies the domain $\mathcal{F} = \mathcal{O} \setminus \overline{\mathcal{S}}$ (see Figure \ref{dessin}).\\
\begin{figure}[!h]
\begin{flushleft}
\hspace*{50pt}\scalebox{0.4}{
\begin{picture}(0,0)%
\includegraphics{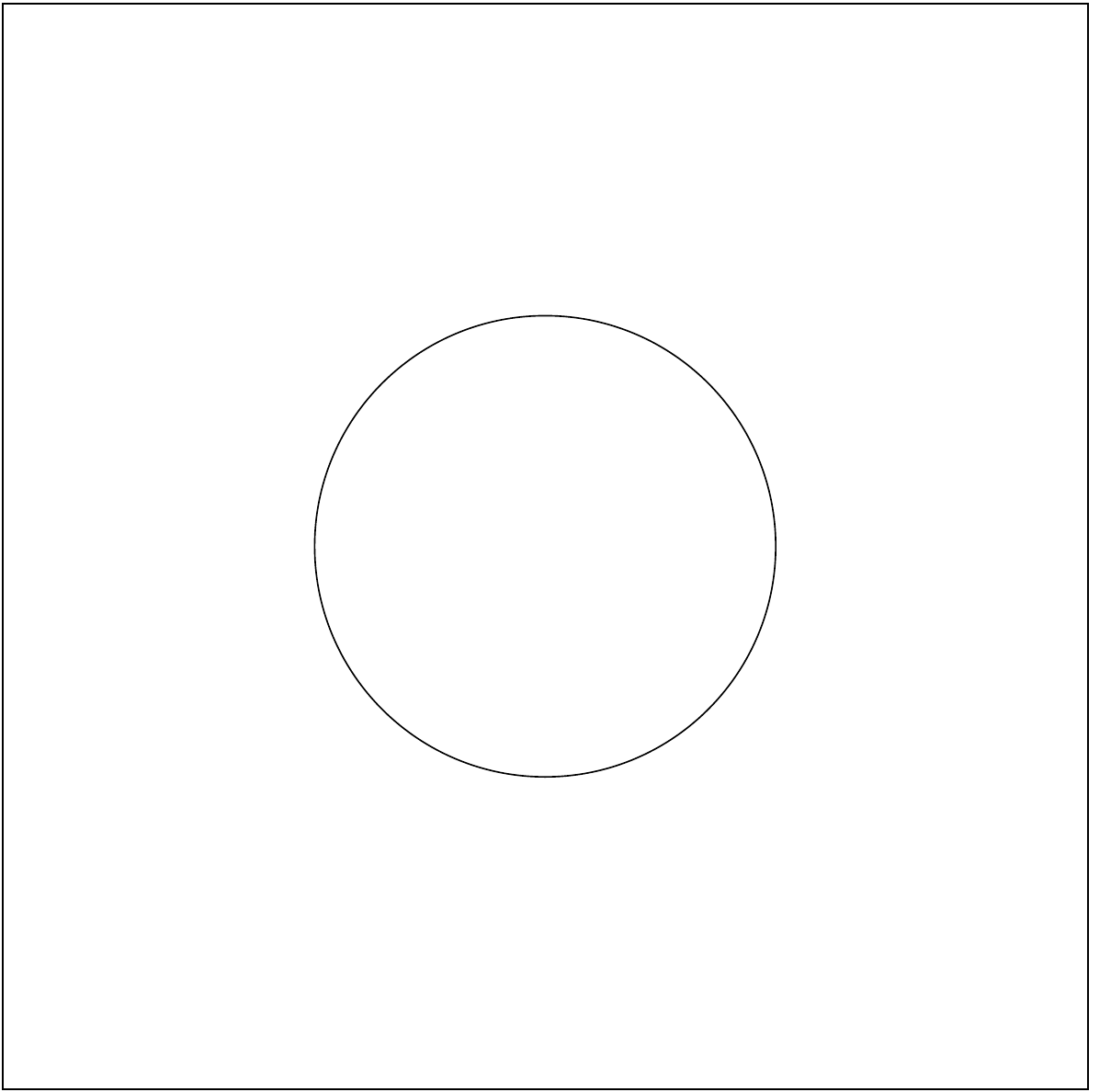}%
\end{picture}%
\setlength{\unitlength}{4144sp}%
\begin{picture}(5517,5424)(889,-5473)
\put(1621,-961){\makebox(0,0)[lb]{\smash{{\SetFigFont{20}{24.0}{\familydefault}{\mddefault}{\updefault}{\color[rgb]{0,0,0}$\mathcal{F}$}%
}}}}
\put(3286,-2716){\makebox(0,0)[lb]{\smash{{\SetFigFont{20}{24.0}{\familydefault}{\mddefault}{\updefault}{\color[rgb]{0,0,0}$\mathcal{S}$}%
}}}}
\put(6391,-556){\makebox(0,0)[lb]{\smash{{\SetFigFont{20}{24.0}{\familydefault}{\mddefault}{\updefault}{\color[rgb]{0,0,0}$\partial \mathcal{O}$}%
}}}}
\put(4726,-3346){\makebox(0,0)[lb]{\smash{{\SetFigFont{20}{24.0}{\familydefault}{\mddefault}{\updefault}{\color[rgb]{0,0,0} $\partial \mathcal{S} = \Gamma$}%
}}}}
\put(4262,-2825){\makebox(0,0)[lb]{\smash{{\SetFigFont{20}{24.0}{\familydefault}{\mddefault}{\updefault}{\color[rgb]{0,0,0} $\longleftarrow$}%
}}}}
\put(4420,-2660){\makebox(0,0)[lb]{\smash{{\SetFigFont{20}{24.0}{\familydefault}{\mddefault}{\updefault}{\color[rgb]{0,0,0} $\bn$}%
}}}}
\end{picture}%
}
\vspace*{-130pt}
\begin{flushright}
\begin{minipage}{8cm}
The velocity of the fluid $\bu$ and its pressure $p$ are assumed to satisfy the following Stokes problem
\begin{eqnarray*}
\left\{\begin{array} {rcl}
- \nu \Delta \bu + \nabla p & = & \bff \quad \text{in } \mathcal{F}, \label{system1} \\
\div \ \bu & = & 0 \quad \text{in }\mathcal{F}, \label{system2} \\
\bu & = & 0 \quad \text{on } \p \mathcal{O}, \label{system3} \\
\bu & = & \bg \quad \text{on } \Gamma, \label{system5}
\end{array} \right.
\end{eqnarray*}
with $\bff\in \mathbf{L}^2(\mathcal{F})$, $\bg\in \mathbf{H}^{1/2}(\Gamma)$, and $\displaystyle \int_{\Gamma} \bg \cdot \bn \d \Gamma =0$.
\end{minipage}
\end{flushright}
\end{flushleft}
\end{figure}
\begin{figure}[!h]
\caption{Domains for the fluid and the solid.\hspace*{230pt}} \label{dessin}
\end{figure}

The Dirichlet condition on the interface $\Gamma$ is nonhomogeneous. The purpose of this work is to impose this condition via a multiplier, and to get a good approximation of the latter for a geometry independent of the mesh. This multiplier represents the normal trace of the Cauchy stress tensor, that is to say
\begin{eqnarray*}
\blambda(\bu,p)  =  \sigma(\bu,p)\bn =  2\nu D(\bu)\bn -p\bn,
\text{ where }
D(\bu)  =  \frac{1}{2} \left(\nabla \bu + \nabla \bu^T \right).
\end{eqnarray*}
The vector $\bn$ denotes the outward unit normal vector to $\p \mathcal{F}$ (see Figure \ref{dessin}). A weak solution of the Stokes system written above can be viewed as a saddle point of the following Lagrangian functional
\begin{eqnarray*}
L_0(\bu,p,\blambda) & = & \nu \int_{\mathcal{F}} \left|D(\bu) \right|^2\d \mathcal{F} - \int_{\mathcal{F}}p\div\ \bu \d \mathcal{F} - \int_{\mathcal{F}}\bff\cdot \bu\d \Gamma - \int_{\Gamma} \blambda \cdot (\bu - \bg) \d \Gamma. \nonumber \\ \label{Lag0}
%\langle \blambda ; \bu-\bg \rangle_{\mathbf{H}^{-1/2}(\Gamma);\mathbf{H}^{1/2}(\Gamma)} . \nonumber \\ \label{Lag0}
\end{eqnarray*}
In order to stabilize the convergence for the multiplier $\blambda$ associated with the Dirichlet condition on $\Gamma$, we carry out an augmented Lagrangian method, like the one initially introduced in \cite{Barbosa1, Barbosa2}. It consists in adding to the Lagrangian functional introduced above a quadratic term, as follows
\begin{eqnarray*}
L(\bu,p,\blambda) & = & L_0(\bu,p,\blambda) - \frac{\gamma}{2}\int_{\Gamma}\left|\blambda - \sigma(\bu,p)\bn \right|^2\d \Gamma \label{Lag1}.
\end{eqnarray*}
The coefficient $\gamma >0$ will be chosen later on. We introduce the following functional spaces
\begin{eqnarray*}
\mathbf{V} & = &\left\{ \bv\in \mathbf{H}^1(\mathcal{F}) \mid \bv=0 \text{ on } \p \mathcal{O} \right\},\\
Q & = & \L^2_0(\mathcal{F}) = \left\{p \in \L^2(\mathcal{F}) \mid \displaystyle \int_{ \mathcal{F}}p\ \d \mathcal{F} = 0 \right\},\\
\mathbf{W} & = & \mathbf{H}^{-1/2}(\Gamma) = \left(\mathbf{H}^{1/2}(\Gamma) \right)'.
\end{eqnarray*}
and we consider the variational problem derived from the Lagrangian $L$ as follows:\footnote{Note that we should assume some additional smoothness in the expression of $L_0$ to make sense, for example $\bu\in \mathbf{H}^2(\mathcal{F})$, $p\in \H^1(\mathcal{F})$, $\blambda\in \mathbf{L}^2(\Gamma)$. The exact solution normally has this smoothness provided that $\bff \in \mathbf{L}^2(\mathcal{F})$ and $\bg \in \mathbf{H}^{3/2}(\Gamma)$.}
\begin{eqnarray*}
& & \text{Find } (\bu,p,\blambda) \in \mathbf{V} \times Q \times \mathbf{W} \text{ such that } \nonumber \\
& & \left\{ \begin{array} {lll}
\mathcal{A}((\bu,p,\blambda);\bv)  = \mathcal{L}(\bv) & \quad & \forall \bv \in \mathbf{V}, \\
\mathcal{B}((\bu,p,\blambda);q)  = 0 & \quad & \forall q \in Q, \\
\mathcal{C}((\bu,p,\blambda);\bmu) = \mathcal{G}(\bmu), & \quad & \forall \bmu \in \mathbf{W},
\end{array} \right. \label{FVaugmented}
\end{eqnarray*}
where
\begin{eqnarray*}
\mathcal{A}((\bu,p,\blambda);\bv) & = & 2\nu\int_{\mathcal{F}}D(\bu):D(\bv)\d \mathcal{F} - \int_{\mathcal{F}}p\div \ \bv\d \mathcal{F}  - \int_{\Gamma} \blambda\cdot \bv\d \Gamma \\
& & -4\nu^2\gamma \int_{\Gamma}\left( D(\bu)\bn \right)\cdot \left( D(\bv)\bn \right)\d \Gamma +2\nu \gamma \int_{\Gamma}p \left(D(\bv)\bn\cdot \bn \right)\d \Gamma +2\nu \gamma \int_{\Gamma} \blambda \cdot \left(D(\bv)\bn\right)\d \Gamma , \\
\mathcal{B}((\bu,p,\blambda);q) & = & - \int_{\mathcal{F}}q\div\ \bu\d \mathcal{F} +2\nu \gamma \int_{\Gamma}q\left(D(\bu)\bn\cdot \bn \right)\d \Gamma -\gamma \int_{\Gamma}pq \d \Gamma - \gamma \int_{\Gamma} q\blambda \cdot \bn\d \Gamma , \\
\mathcal{C}((\bu,p,\blambda);\bmu) & = & -\int_{\Gamma} \bmu \cdot \bu \d \Gamma +2\nu \gamma \int_{\Gamma}\bmu \cdot (D(\bu)\bn)\d \Gamma -\gamma \int_{\Gamma}p(\bmu\cdot \bn)\d \Gamma - \gamma \int_{\Gamma} \blambda \cdot \bmu \d \Gamma , \\
\mathcal{L}(\bv) & = & \int_{\mathcal{F}}\bff\cdot \bv\d \mathcal{F}, \label{defaL} \qquad \quad
\mathcal{G}(\bmu) \ \ = \ \ -\int_{\Gamma} \bmu \cdot \bg \d \Gamma. \label{defG}
\end{eqnarray*}

\section{Description of the method} \label{sec2}

\subsection{Fictitious domain approach}
\begin{figure}[!h]
\begin{center}
\hspace*{2cm} \includegraphics[trim = 9cm 2cm 2cm 1cm, clip, scale= 0.50]{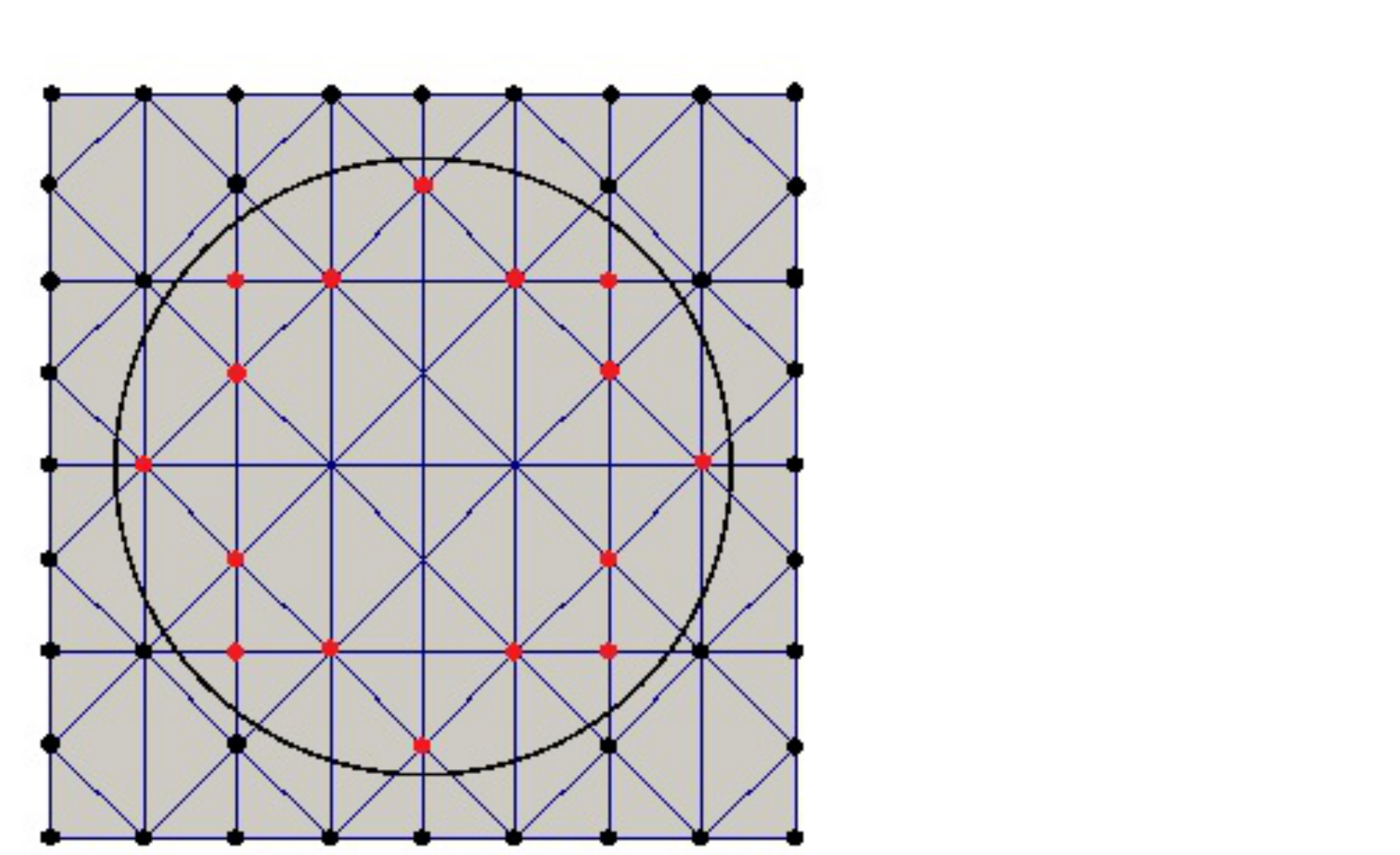}\hspace*{-2cm}  
\includegraphics[trim = 9cm 2cm 2cm 1cm, clip, scale= 0.50]{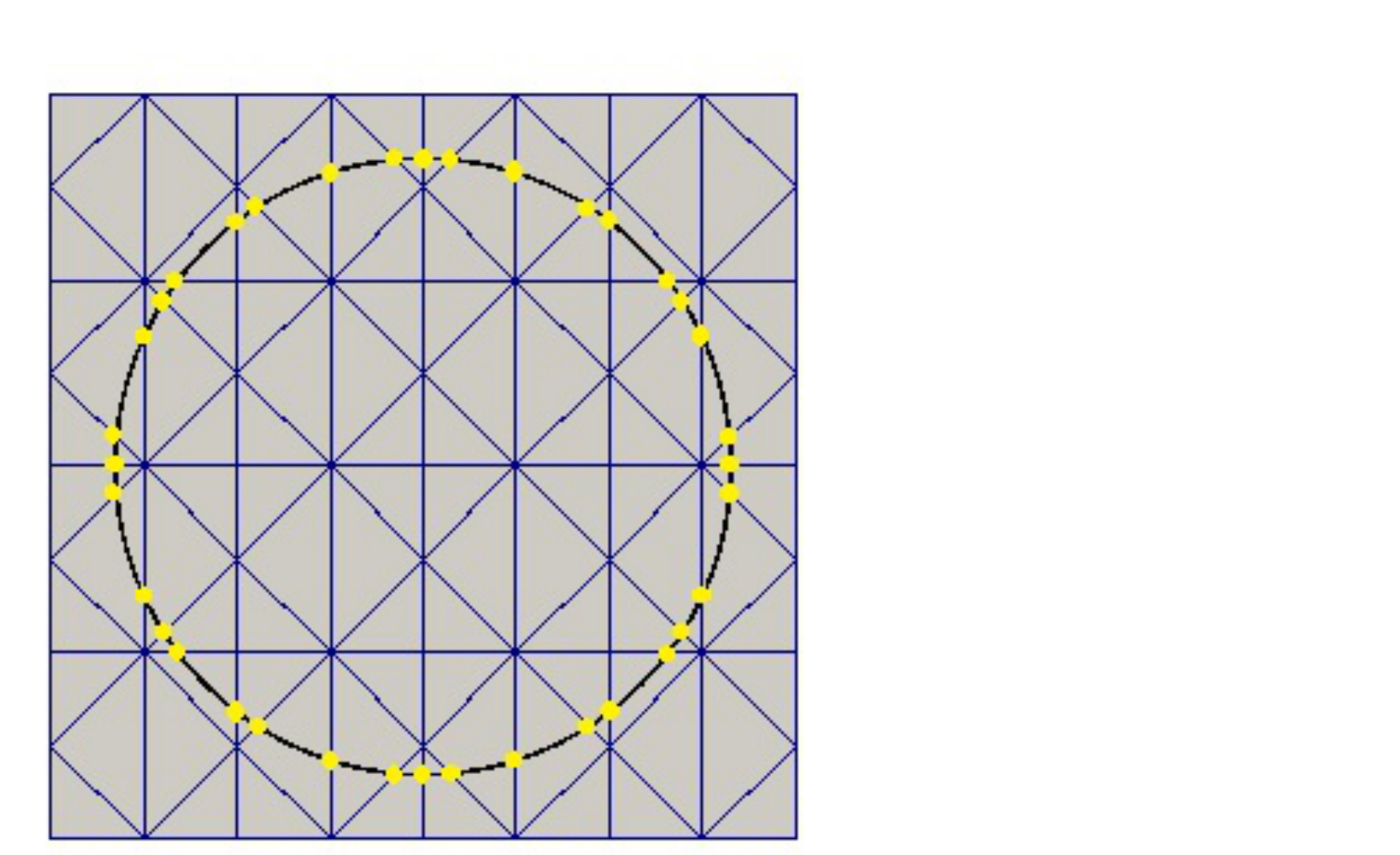}
\caption{Mesh on a fictitious domain. In the first figure: Standard degrees of freedom (black), virtual ones (red), remaining ones are removed. In the second figure: Bases nodes used for the multiplier space (yellow).} \label{mesh}
\end{center}
\end{figure}

The fictitious domain for the fluid is considered on the whole domain  $\mathcal{O}$. Let us introduce three discrete finite element spaces, {$\tilde{\mathbf{V}}^h \subset \mathbf{H}^1(\mathcal{O})$, $\tilde{Q}^h \subset \L^2_0(\mathcal{O})$ and $\tilde{\mathbf{W}}^h \subset \mathbf{L}^2(\mathcal{O})$}. Since $\mathcal{O}$ can be a rectangular domain, \textcolor{black}{these} spaces can be defined on the same structured mesh, that can be chosen uniform (see Figure \ref{mesh}). The construction of the mesh is highly simplified (no particular mesh is required). We set
\begin{eqnarray}
\tilde{\mathbf{V}}^h & = & \left\{\bv^h \in C(\overline{\mathcal{O}})\mid \bv^h_{\left| \p \mathcal{O}\right.} = 0, \  \bv^h_{\left| T\right.} \in P(T), \ \forall T \in \mathcal{T}^h \right\}, \label{defvtilde}
\end{eqnarray}
where $P(T)$ is a finite dimensional space of regular functions such that $P(T) \supseteq P_k(T)$ for some integer $k \geq 1$. For more details, see \cite{Ern} for instance. The mesh parameter stands for $\displaystyle h = \max_{T\in \mathcal{T}^h} h_T$, where $h_T$ is the diameter of $T$. We define
\begin{eqnarray*}
\mathbf{V}^h := \tilde{\mathbf{V}}^h_{\left| \mathcal{F} \right.}, \quad Q^h := \tilde{Q}^h_{\left|\mathcal{F}\right.}, \quad  \mathbf{W}^h := \tilde{\mathbf{W}}^h_{\left| \Gamma \right.},
\end{eqnarray*}
which are natural discretizations of $\mathbf{V}$, $\L^2(\mathcal{F})$ and $\mathbf{H}^{-1/2}(\Gamma)$ respectively. 
This approach is equivalent to eXtended Finite Element Method as proposed in \cite{Choi2010} or \cite{Gerstenberger2008} where the standard Finite Element Method basis functions are multiplied by the Heaviside function ($H({\bf x}) = 1$ for ${\bf x} \in \mathcal{F}$ and $H({\bf x})=0$ for ${\bf x}\in \mathcal{O} \setminus \mathcal{F}$) and the products are substituted in the variational formulation of the problem. Thus the degrees of freedom inside the fluid domain $\mathcal{F}$ are used in the same way as in the standard Finite Element Method, whereas the degrees of freedom in the solid domain $\mathcal{S}$ at the vertices of the elements cut by the interface (the so called virtual degrees of freedom) do not define the field variable at these nodes, but they are necessary to define the fields on $\mathcal{F}$ and to compute the integrals over $\mathcal{F}$. The remaining degrees of freedom, corresponding to the basis functions with support completely outside of the fluid, are eliminated (see Figure 2). We refer to the papers mentioned above for more details.

\subsection{The discrete stabilized formulation}
Let us choose $ \gamma = \gamma_0 h$, where the constant $\gamma_0 >0$ has to be chosen sufficiently large for stabilizing, but not too much in order to keep the system coercive. The discrete problem can be rewritten in the following compact form:
\begin{eqnarray*}
& & \text{Find $(\bu^h,p^h,\blambda^h) \in \mathbf{V}^h \times Q^h \times \mathbf{W}^h$ such that} \\
& & \mathcal{M}((\bu^h,p^h,\blambda^h);(\bv^h,q^h,\bmu^h))  =
{\cal H}(\bv^h,q^h,\bmu^h), \quad \forall (\bv^h,q^h,\bmu^h) \in \mathbf{V}^h \times Q^h \times \mathbf{W}^h,
\end{eqnarray*}
where
\begin{eqnarray*}
\mathcal{M}((\bu,p,\blambda);(\bv,q,\bmu)) &=&2\nu \int_{\mathcal{F}}D(\bu):D(\bv)\d \mathcal{F}-\int_{\mathcal{F}}(p\div \ \bv+q\div \ \bu) \d \mathcal{F}-\int_{\Gamma }(\blambda\cdot \bv+\bmu\cdot \bu)\d\Gamma  \\
&&-\gamma_0 h \int_{\Gamma }(2\nu D(\bu)\bn-p\bn-\blambda)\cdot \left( 2\nu D(\bv%
)\bn-q\bn-\bmu\right) \d\Gamma,\\
\mathcal{H}(\bv,q,\bmu) & = & \int_{\mathcal{F}} \bff\cdot \bv \d \Gamma  -\int_{\Gamma} \bmu \cdot \bg \d \Gamma.
\end{eqnarray*}

\subsection{Convergence}
Let us give some assumptions we need for stating an important theoretical result \textcolor{black}{(See \cite{Court}, section 4.2)}:
\begin{description}
\item[A1] For all $\bv^{h}\in \mathbf{V}^{h}$ one has\footnote{\textcolor{black}{We denote by $C$ a generic positive constant which does not depend on the mesh size $h$.}}
\begin{eqnarray*}
h\|D(\bv^{h})n\|_{\mathbf{L}^2(\Gamma)}^{2} & \leq & C\|\bv^{h}\|_{\mathbf{V}}^{2}.
\end{eqnarray*}
\item[A2] For all $q^{h}\in Q^{h}$ one has%
\begin{eqnarray*}
h\|q^{h}\|_{\L^2(\Gamma)}^{2} & \leq & C\|q^{h}\|_{\L^2(\mathcal{F})}^{2}.
\end{eqnarray*}
\item[A3] One has the following inf-sup condition for the velocity-pressure pair of finite
element spaces
\begin{eqnarray*}
\inf_{q^{h}\in Q^{h}}\sup_{\bv^{h}\in \mathbf{V}_{0}^{h}}\frac{b(\bv^h,q^h)}{\|q^{h}\|_{\L^2(\mathcal{F})}\|\bv^{h}\|_{\mathbf{V}}} & \geq & \text{\textcolor{black}{$C$}.}
\end{eqnarray*}
%with $\beta >0$ independent of $h$.
\end{description}

\begin{lemma} \label{lemmainfsup}
Under assumptions {\bf A1}--{\bf A3} given above, there exists for $\gamma_0$ small enough a mesh-independent constant $c>0$ such that
\begin{eqnarray*}
\inf_{(\bu^{h},p^{h},\blambda^{h})\in \mathbf{V}^{h}\times Q^{h}\times \mathbf{W}^{h}} \sup_{(\bv^{h},q^{h},\bmu^{h})\in \mathbf{V}^{h}\times Q^{h}\times \mathbf{W}^{h}}\frac{%
\mathcal{M}((\bu^{h},p^{h},\blambda^{h});(\bv^{h},q^{h},\bmu ^{h}))}{%
|||\bu^{h},p^{h},\blambda^{h}|||\,|||\bv^{h},q^{h},\bmu ^{h}|||} & \geq & \text{\textcolor{black}{$C$}},
\end{eqnarray*}%
where the triple norm is defined by%
\begin{eqnarray*}
|||\bu,p,\blambda |||=\left( \|\bu\|_{\mathbf{V}}^{2}+\|p\|_{\L^2(\mathcal{F})}^{2}+h\|D(\bu)\bn%
\|_{\mathbf{L}^2(\Gamma) }^{2}+h\|p\|_{\L^2(\Gamma) }^{2}+h\|\blambda \|_{\mathbf{L}^2(\Gamma) }^{2}+%
\frac{1}{h}\|\bu\|_{\mathbf{L}^2(\Gamma) }^{2}\right) ^{1/2}.
\end{eqnarray*}%
\end{lemma}

A consequence of this lemma is \textcolor{black}{the} following abstract error estimate
\begin{eqnarray*}\label{eststab}
&&\max ( \|\bu-\bu^{h}\|_{\mathbf{V}} , \|p-p^h\|_{\L^2(\mathcal{F})}, h\|\blambda -\blambda^h\|_{\mathbf{L}^2(\Gamma)} )
 \le
|||\bu-\bu^{h},p-p^{h},\blambda-\blambda ^{h}|||
\\
&&\qquad \le C\left(
h^{k_u}\|\bu\|_{\mathbf{H}^{k_u+1}(\mathcal{F})} + h^{k_p+1}\|p\|_{\H^{k_p+1}(\mathcal{F})} +
h^{k_\lambda+1}\|\blambda\|_{\mathbf{H}^{k_{\lambda}+1/2}(\Gamma)}
\right),
\notag
\end{eqnarray*}
where $k_u$, $k_p$ and $k_\lambda$ are the degrees of finite elements used for velocity, pressure and multiplier $\blambda$ respectively.

\section{Robustness with respect to the geometry} \label{sec3}
In a framework where the solid will be led to move in the fluid domain, we need to consider different types of intersection between the level-set and the regular mesh, and then examining the behavior of our method. For instance we can compute the $\mathbf{L}^2(\Gamma)$ relative errors on the multiplier $\blambda$, for different positions of the center of the solid, with and without the stabilization technique. The perspective is to anticipate the behavior of the method in an unsteady case, and these tests enable us to avoid to handle in a first time the complexity of a full unsteady problem.\\
For $h = 0.05$ and the finite elements triplet P2/P1/P0, we consider the solid as a circle, and we make the abscissa of the center of the circle - denoted by $\x_C$ - vary between 0.5 and 0.7 (with a step equal to 0.0005) in a box $[0,1]\times [0,1]$. The variations of the relative error (in \%) on $\blambda$ are represented in blue (without stabilization) and in red (with stabilization).

\begin{minipage}{21cm}
\begin{center}
\hspace*{-5cm}\includegraphics[trim = 0cm 7cm 0cm 7cm, clip, scale=0.6]{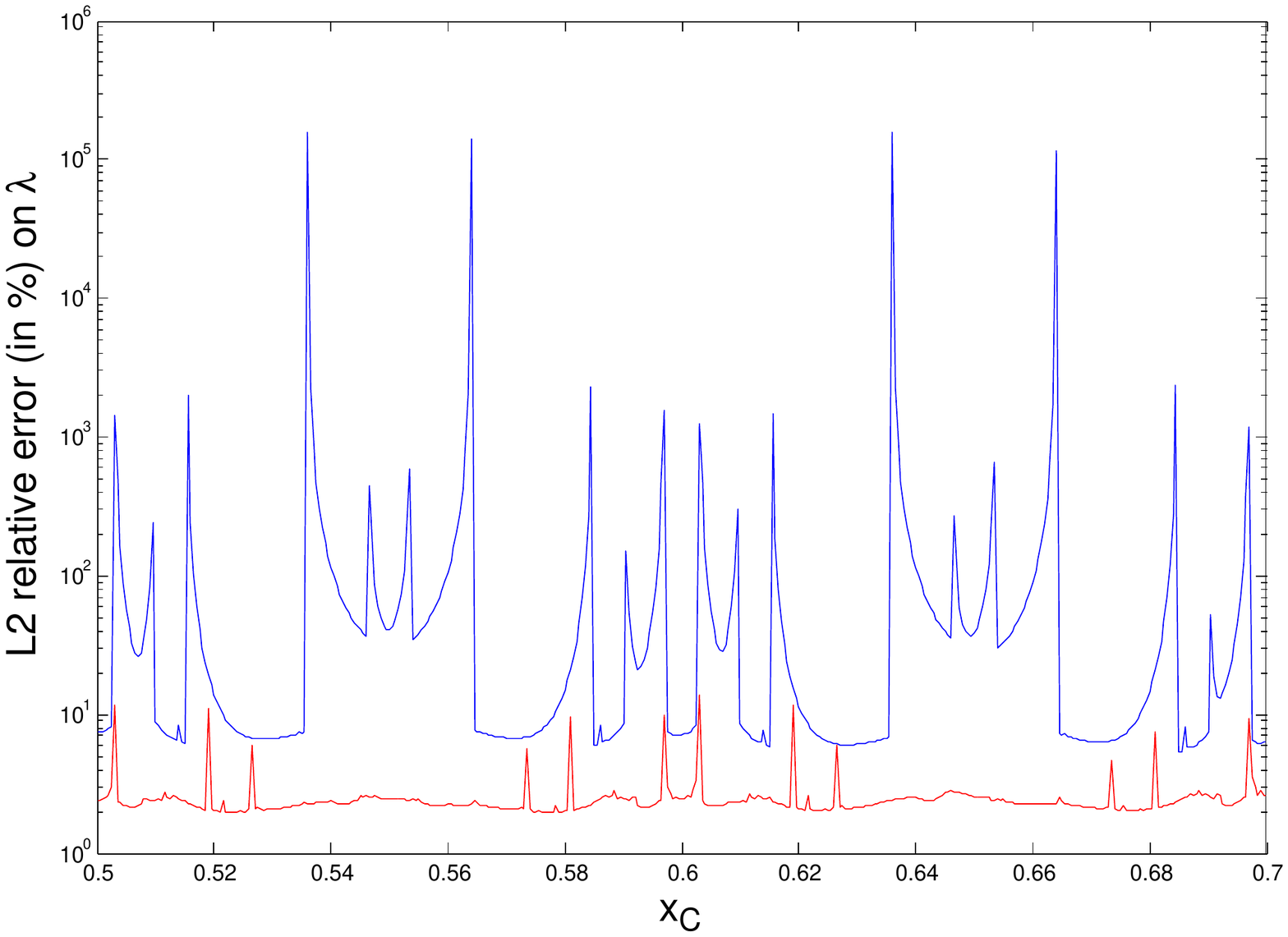}
\end{center}
\end{minipage}

\vspace*{-0.5cm}

\textcolor{black}{
\begin{figure}[!h]
\centering\caption{Behavior of the $\mathbf{L}^2(\Gamma)$ relative error on $\blambda$ (in semi-log scale), in red with the stabilization technique (with $\gamma_0 = 0.05$), in blue without.}
%, and $\max(C_{\bu},C_{p})$ (in green), in function of $x_C$.}
\label{moverror}
\end{figure}
}
\FloatBarrier

In these tests the relevance of our approach using the stabilization technique is underlined when the intersection between the level-set and the mesh varies. Without stabilization the errors are huge in many cases (see the curve in blue), whereas the robustness of the stabilization technique is demonstrated with regards to the constancy of the relative errors (see the curve in red).

\section{Application to a fluid-structure interaction problem:\\ Free fall of an ellipse} \label{sec4}
In order to illustrate our main purpose, that is to say performing numerical tests of fluid-structure problems for time-depending geometries, we propose to simulate the free fall of an ellipsoidal solid, in a viscous incompressible fluid, submitted to its own weight, to the forces that the surrounding fluid applies on its surface (represented by the quantity $\blambda = \sigma(\bu,p)\bn$), and to the Archimedes force. 

\subsection{Description of the model}
While moving, the solid occupies a time-depending time domain $\mathcal{S}(t)$. The remaining domain $\mathcal{F}(t) = \mathcal{O} \setminus \overline{\mathcal{S}(t)}$ is the one occupied by the fluid flow. The displacement of the rigid solid can be represented by a Lagrangian mapping $X_{S}$ given by
\begin{eqnarray*}
\mathcal{S}(t) & = & X_{S}(\mathcal{S}(0)), \\
X_S(\y,t) & = & \bh(t) + \mathbf{R}(t)(\y-\bh(0)), \quad \y \in \mathcal{S}(0),
\end{eqnarray*}
where $\bh(t)$ is the position of the gravity center of the solid, and $\mathbf{R}(t)$ is a rotation written is 2D as $\displaystyle \left( \begin{matrix} c & -s \\ s & c \end{matrix}\right)$, with $c=\cos(\theta(t))$ and $s=\sin(\theta(t))$. The Reynolds number we choose to consider here is intermediate, which means that the fluid has to obey the incompressible Navier-Stokes system
\begin{eqnarray*}
\frac{\p \bu}{\p t} - \nu \Delta \bu + (\bu\cdot \nabla )\bu + \nabla p  =  0 , & \quad & \x \in \mathcal{F}(t),\quad t\in (0,T), \label{introfr_sys1} \\
\div \ \bu   =  0 ,  & \quad & \x\in \mathcal{F}(t), \quad t\in (0,T), \label{introfr_sys2}
\end{eqnarray*}
to which we add Dirichlet conditions, in particular the one imposed by the equality of velocities at the fluid-solid interface $\p \mathcal{S}(t)$:
\begin{eqnarray*}
\bu  =  0 , & \quad & \x \in \p \mathcal{O} ,\quad t\in (0,T),  \\
\bu  =  \bh'(t) +  \bomega (t)\wedge(\x-\bh(t))  , & \quad &  \x\in \p \mathcal{S}(t),\quad t\in(0,T).
\end{eqnarray*}
\textcolor{black}{The quantity $\bomega = \theta'$ denotes the angular velocity of the solid.} Note that the quantities $\bh$ and $\theta$ are also unknowns of the problem. Their dynamics is coupled with the fluid forces by the Newton laws given by
\begin{eqnarray*}
m \bh''(t)  =  - \int_{\p \mathcal{S}(t)} \sigma(\bu,p) \bn \d \Gamma -m\bg +m_{a}\bg
  , & \quad & t\in (0,T),  \\
\left(I\bomega\right)' (t)  = - \int_{\p \mathcal{S}(t)} (\x-h(t))\wedge  \sigma(\bu,p) \bn \d \Gamma, & \quad & t\in(0,T),
\end{eqnarray*}
where $m$, $I$ denote the mass and the inertia moment of the solid respectively, $\bg$ is the gravity field and $m_{a}$ comes from the Archimedes' principle corresponding to the fluid mass displaced by the volume of the solid.\\
The initial velocities are chosen to be equal to zero, as well for the fluid as for the solid. The latter is dropped in a box $[0,1]\times [0,2.5]$ while being inclined of an angle $\theta_0$ equal to $1.40$ rad. (see Figure \ref{panorama}). The dimensions of its semi-major axis and semi-minor axis are respectively
\begin{eqnarray*}
a = 0.24, & \quad & b = 0.08.
\end{eqnarray*}
The viscosity of the fluid and the mass of the solid are respectively chosen as
\begin{eqnarray*}
\nu = 1.00, & \quad & m = 20.0.
\end{eqnarray*}
The mesh on the whole domain is triangular and based on 40 subdivisions in horizontal direction and 100 subdivisions in vertical direction. For the triplet $(\bu, p, \blambda)$, we choose the classical finite element P2/P1/P0 and the stabilization parameter $\gamma = h\times \gamma_0$ for $\gamma_0 = 0.05$.

\subsection{The discrete problem}
If the degrees of freedom of $\bu^h$, $p^h$ and $\blambda^h$ are represented by the vectors $\boldsymbol{U}$, $\boldsymbol{P}$ and $\boldsymbol{\Lambda}$ respectively, the discrete problem can be formulated as follows

\begin{eqnarray*}
M_{{\bf u}{\bf u}} \frac{\mbox{d} \boldsymbol{U(t)}}{\mbox{d} t} +A_{{\bf u}{\bf u}} \boldsymbol{U}(t) + N( \boldsymbol{U}(t)) \boldsymbol{U}(t) +  A_{{\bf u}p}\boldsymbol{P}(t) +  A_{{\bf u}{\bf \lambda}}\boldsymbol{\Lambda}(t)  = 0, \\
A^T_{{\bf u}p}  \boldsymbol{U}(t) +  A_{pp} \boldsymbol{P}(t)  + A_{p{\bf \lambda}}\boldsymbol{\Lambda}(t) = 0,  \\
A^T_{{\bf u}{\bf \lambda}}  \boldsymbol{U}(t)  +  A^T_{p{\bf \lambda}} \boldsymbol{P}(t)  +  A_{{\bf \lambda} {\bf \lambda}}\boldsymbol{\Lambda}(t) = \boldsymbol{G},\\
\displaystyle m{\bf h}''(t) = M_{ \bf \lambda} \boldsymbol{\Lambda}(t) -m {\bf g} + m_a {\bf g},&&\\
\displaystyle \text{{$I\theta''(t) =  M_{ \bf \lambda}
\left[ ({\bf x}-{\bf h}(t))^{\bot} \cdot \boldsymbol{\Lambda}(t)\right]$}},	&&
\end{eqnarray*}
where the matrices which appear above are discretizations of the following mappings
\begin{eqnarray*}
\mathcal{M}_{{\bf u}{\bf u}} : ({\bf u},{\bf v}) & \longmapsto & \int_{\mathcal{F}}{\bf u}.{\bf v}\mathrm{d}\mathcal{F}, \qquad
\mathcal{M}_{\blambda} : {\blambda}   \longmapsto  -\int_{\partial \mathcal{S}(t)}{\blambda} \mathrm{d}\Gamma,\\
\mathcal{A}_{{\bf u}{\bf u}} : ({\bf u},{\bf v}) & \longmapsto & 2\nu\int_{\mathcal{F}}D({\bf u}):D({\bf v})\mathrm{d}\mathcal{F} - 4\nu^2\gamma \int_{\Gamma}\left( D({\bf u})\bn \right)\cdot \left( D({\bf v})\bn\right)\mathrm{d}\Gamma , \\
\mathcal{A}_{{\bf u}p} : ({\bf v},p) & \longmapsto & - \int_{\mathcal{F}}p\div \ {\bf v} \mathrm{d}\mathcal{F} + 2\nu \gamma \int_{\Gamma}p \left(D({\bf v})\bn\cdot \bn \right)\mathrm{d}\Gamma , \\
\mathcal{A}_{{\bf u}{\blambda}} : ({\bf u},{\blambda}) & \longmapsto & - \int_{\Gamma} {\blambda}\cdot {\bf v}\mathrm{d}\Gamma + 2\nu \gamma \int_{\Gamma} {\blambda} \cdot \left(D({\bf v})\bn\right)\mathrm{d}\Gamma , \\
\mathcal{A}_{pp} : (p,q) & \longmapsto & -\gamma \int_{\Gamma}pq \mathrm{d}\Gamma , \hspace*{0.5cm} 
\mathcal{A}_{p{\blambda}} : (q,{\blambda})  \longmapsto  - \gamma \int_{\Gamma} q{\blambda} \cdot \bn\mathrm{d}\Gamma , \\
\mathcal{A}_{{\blambda}{\blambda}} : ({\blambda},{\bmu}) & \longmapsto & -\gamma \int_{\Gamma} {\blambda} \cdot {\bmu} \mathrm{d}\Gamma, \qquad 
%\mathcal{G} : {\bf \mu}  \longmapsto  - \int_{\Gamma} {\bf \mu} \cdot  (h'(t) + \theta'(t) ({\bf x}-h(t))^{\bot})\mathrm{d}\Gamma, 
\end{eqnarray*}
and the vector  $\boldsymbol{G}$ is the discretization of  $\mathcal{G}$. The term $N(\boldsymbol{U}(t))\boldsymbol{U}(t)$ is a matrix depending on the velocity corresponding to the nonlinear convective term  $\displaystyle \int_{\mathcal{F}} [({\bf u} \cdot \nabla){\bf u}] \cdot{\bf v} \mathrm{d}\mathcal{F}$.

\subsection{Time discretization and treatment of the nonlinearity} \label{time}
We consider an implicit time discretization based on the backward Euler method. We denote by $\boldsymbol{U}^{n+1}$ the solution at the time level $t^{n+1}$ and $dt = t^{n+1}-t^n$ is the time step. Particular attention must be done for a moving particle problem. Indeed, at the time level $t^{n+1}$ the solid occupies the domain $\mathcal{S}(t^{n+1})$ which is different of the previous time level $t^n$. Thus the field variable at time $t^{n+1}$ can become undefined near the interface where there was no fluid flow at time $t^n$ ($\mathcal{S}(t^{n+1})\neq \mathcal{S}(t^{n})$ for the solid and $\mathcal{F}(t^{n+1})\neq \mathcal{F}(t^{n})$ for the fluid). Some {degrees} of freedom (inside the solid) which are not considered at time $t^n$ have to be taken into account at time $t^{n+1}$. In particular the velocity field must be known in such nodes. In this work we impose the velocity to be equal to the motion of the solid.\\
Let us give the algorithm which enables us to compute at the time level $t^{n+1}$ the solution on $\mathcal{F}(t^{n+1})$ represented by $(\boldsymbol{U}^{n+1}, \boldsymbol{P}^{n+1}, \boldsymbol{\Lambda}^{n+1},{\bf h}'^{n+1},{\bf h}^{n+1},\theta'^{n+1},\theta^{n+1})$. Note that at the time level $t^n$ we have access to ($\boldsymbol{U}^n, \boldsymbol{P}^n, \boldsymbol{\Lambda}^n,{\bf h}'^n,{\bf h}^n,\theta'^n,\theta^{n}$) on $\mathcal{F}(t^{n})$.
\begin{small}
\begin{itemize}
\item[1--]  We compute $({\bf h}'^{n+1},\theta'^{n+1})$ such that
\begin{eqnarray*}
\displaystyle m \frac{{\bf h}'^{n+1}-{\bf h}'^n}{dt}  = M_{ \bf \lambda} \boldsymbol{\Lambda}^{n} -m {\bf g} + m_a {\bf g}, & &
\displaystyle I\frac{\theta'^{n+1} - \theta'^n}{dt} = M_{ \bf \lambda}  \left [ ({\bf x}-{\bf h}^{n})^{\bot}\cdot  \boldsymbol{\Lambda}^{n}\right ].
\end{eqnarray*}
\item[2--]  We complete the velocity $\boldsymbol{U}^n$ defined on $\mathcal{F}(t^{n})$ to the full domain by imposing the velocity on each node of $\mathcal{S}(t^{n})$ 
equal to ${\bf h}'^{n+1} + \theta'^{n+1}({\bf x}-{\bf h}^n)^{\bot}$.
\item[3--]  We update the geometry to determine $\mathcal{F}(t^{n+1})$ by computing $({\bf h}^{n+1},\theta^{n+1})$ such that
\begin{eqnarray*}
\displaystyle m \frac{{\bf h}^{n+1}-2{\bf h}^n+{\bf h}^{n-1}}{dt^2}  =  M_{ \bf \lambda} \boldsymbol{\Lambda}^{n} -m {\bf g} + m_a {\bf g}, & &
\displaystyle I\frac{\theta^{n+1} - 2 \theta^n +  \theta^{n-1}}{dt^2} = M_{ \bf \lambda}  \left [ ({\bf x}-{\bf h}^{n})^{\bot}\cdot  \boldsymbol{\Lambda}^{n}\right ].
\end{eqnarray*}
\item[4--] We compute the Dirichlet condition for the velocity at the new interface $\Gamma = \partial \mathcal{S}(t^{n+1})$. 
So we determine $\boldsymbol{G}^{n+1}$ from ${\bf u}_{\Gamma}^{n+1} = {\bf h}'^{n+1} + \theta'^{n+1} ({\bf x}-{\bf h}^{n+1})^{\bot}$.
\item[5--] Finally, we find $(\boldsymbol{U}^{n+1},\boldsymbol{P}^{n+1},\boldsymbol{\Lambda}^{n+1})$ such that
\begin{eqnarray*}
M_{{\bf u}{\bf u}} \frac{\boldsymbol{U}^{n+1} -\boldsymbol{U}^{n} }{dt} +A_{{\bf u}{\bf u}} \boldsymbol{U}^{n+1} + N(\boldsymbol{U}^{n+1}) \boldsymbol{U}^{n+1} +  A_{{\bf u}p}\boldsymbol{P}^{n+1}
+  A_{{\bf u}{\bf \lambda}}\boldsymbol{\Lambda}^{n+1} = 0, \\
A^T_{{\bf u}p}  \boldsymbol{U}^{n+1} +  A_{pp} \boldsymbol{P}^{n+1}  + A_{p{\bf \lambda}}\boldsymbol{\Lambda}^{n+1} = 0,  \\
A^T_{{\bf u}{\bf \lambda}}  \boldsymbol{U}^{n+1}  +  A^T_{p{\bf \lambda}} \boldsymbol{P}^{n+1}  +  A_{{\bf \lambda} {\bf \lambda}}\boldsymbol{\Lambda}^{n+1} = \boldsymbol{G}^{n+1}.
\end{eqnarray*}
At this stage, the solution of the resulting nonlinear algebraic system is achieved by a Newton technique.
\end{itemize}
\end{small}

\subsection{Simulation}
For a time step chosen to be equal to $10^{-3}$, the amplitude of the velocity and the evolution of the solid (its position and its orientation) is represented in Figure \ref{panorama}.
\begin{center}
\begin{tabular} {ccccc}
\includegraphics[trim = 16.5cm 1cm 16.5cm 1cm, clip, scale=0.25]{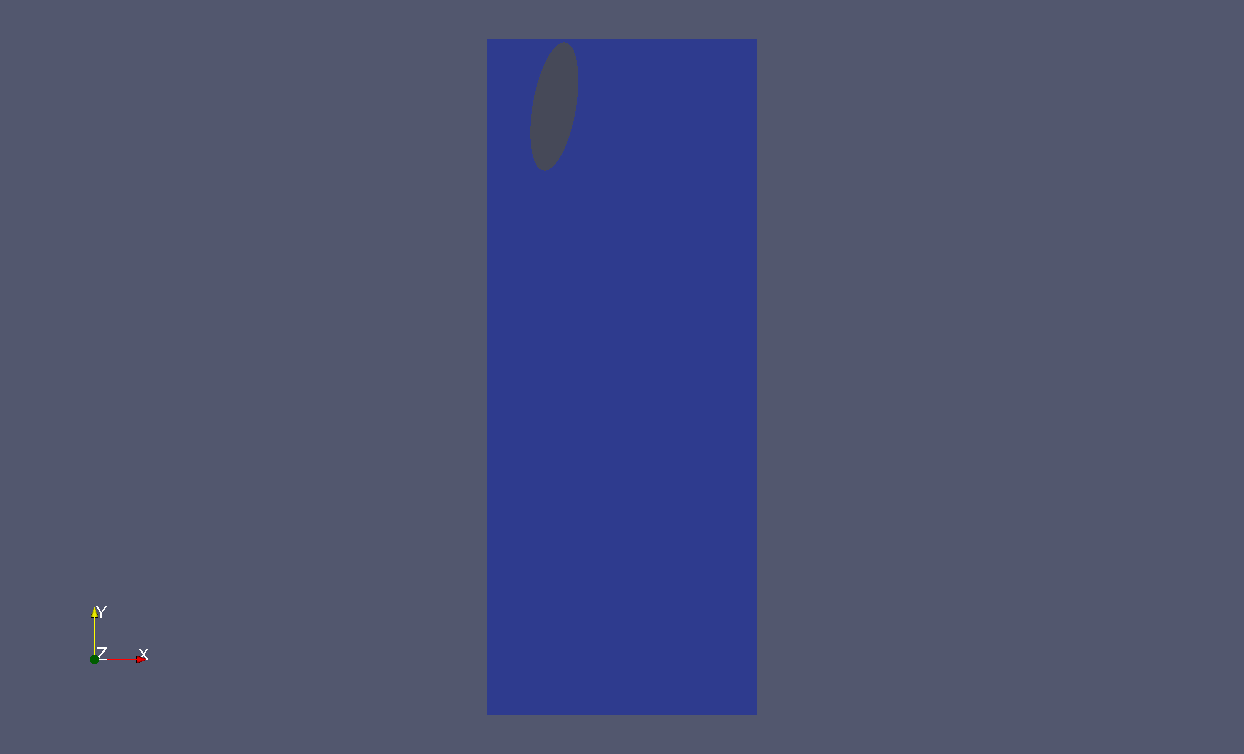} &
\includegraphics[trim = 16.5cm 1cm 16.5cm 1cm, clip, scale=0.25]{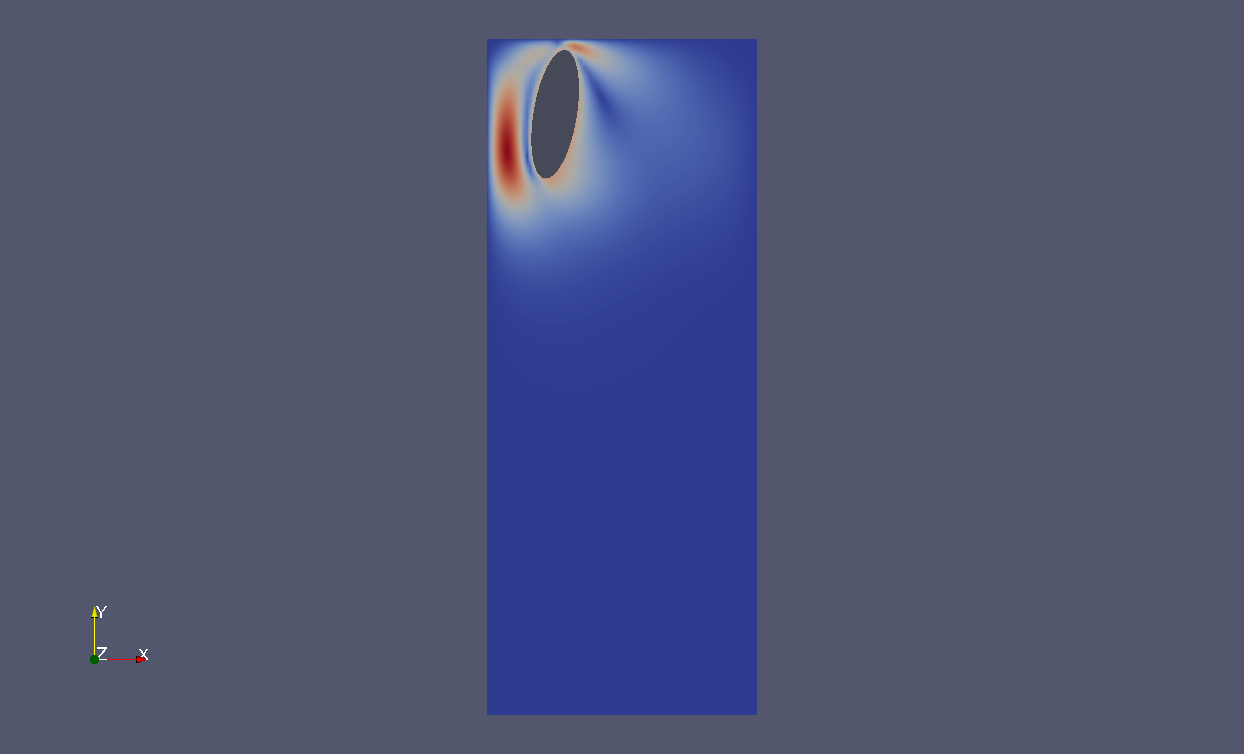} &
\includegraphics[trim = 16.5cm 1cm 16.5cm 1cm, clip, scale=0.25]{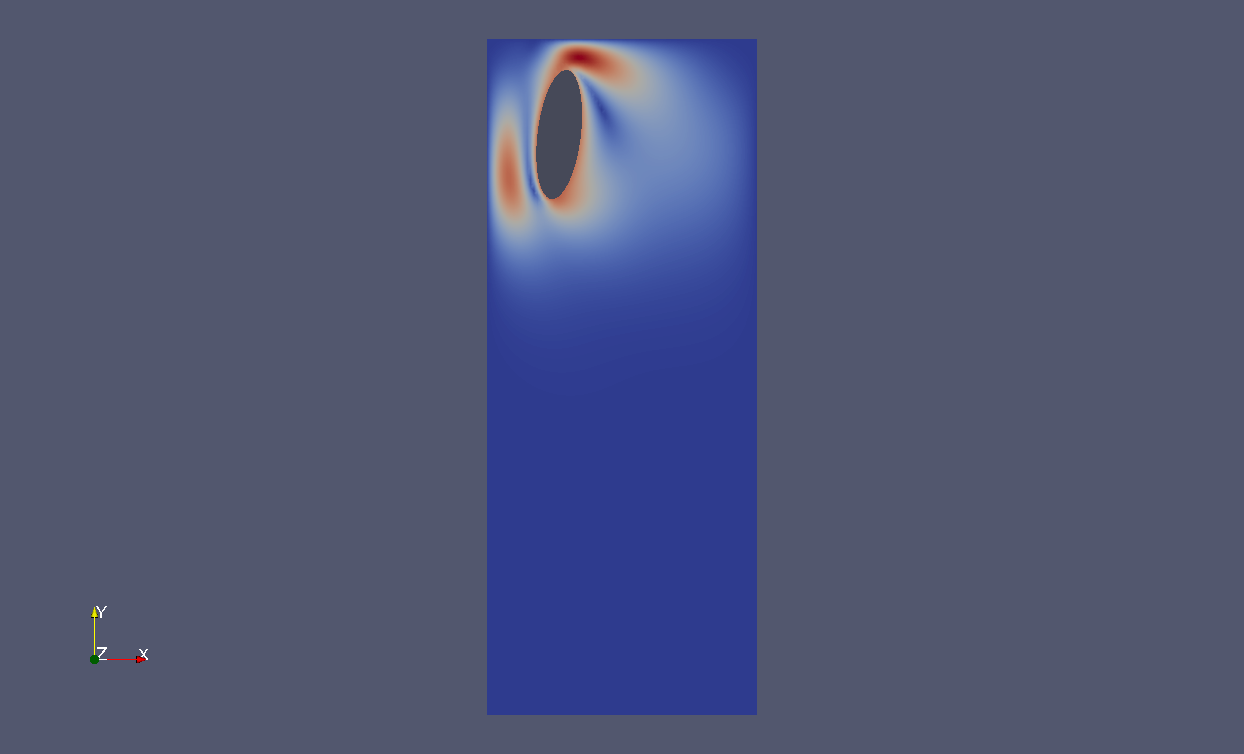} &
\includegraphics[trim = 16.5cm 1cm 16.5cm 1cm, clip, scale=0.25]{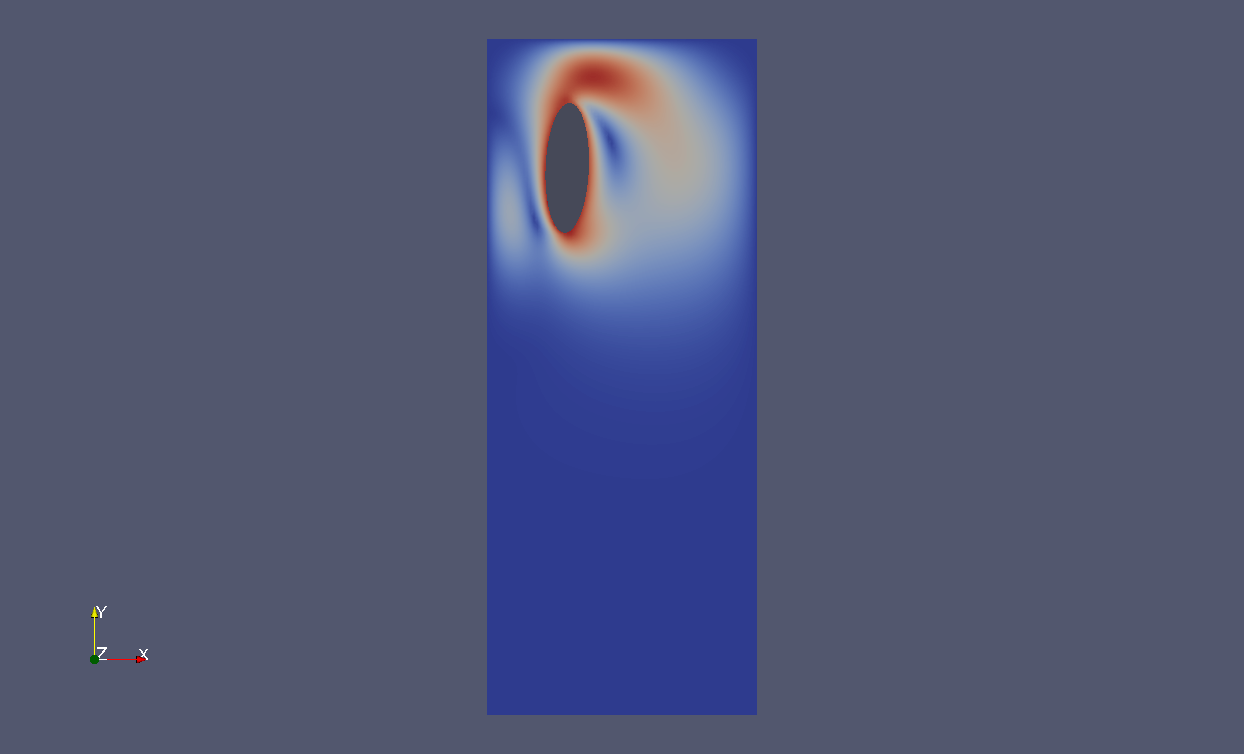} &
\includegraphics[trim = 16.5cm 1cm 16.5cm 1cm, clip, scale=0.25]{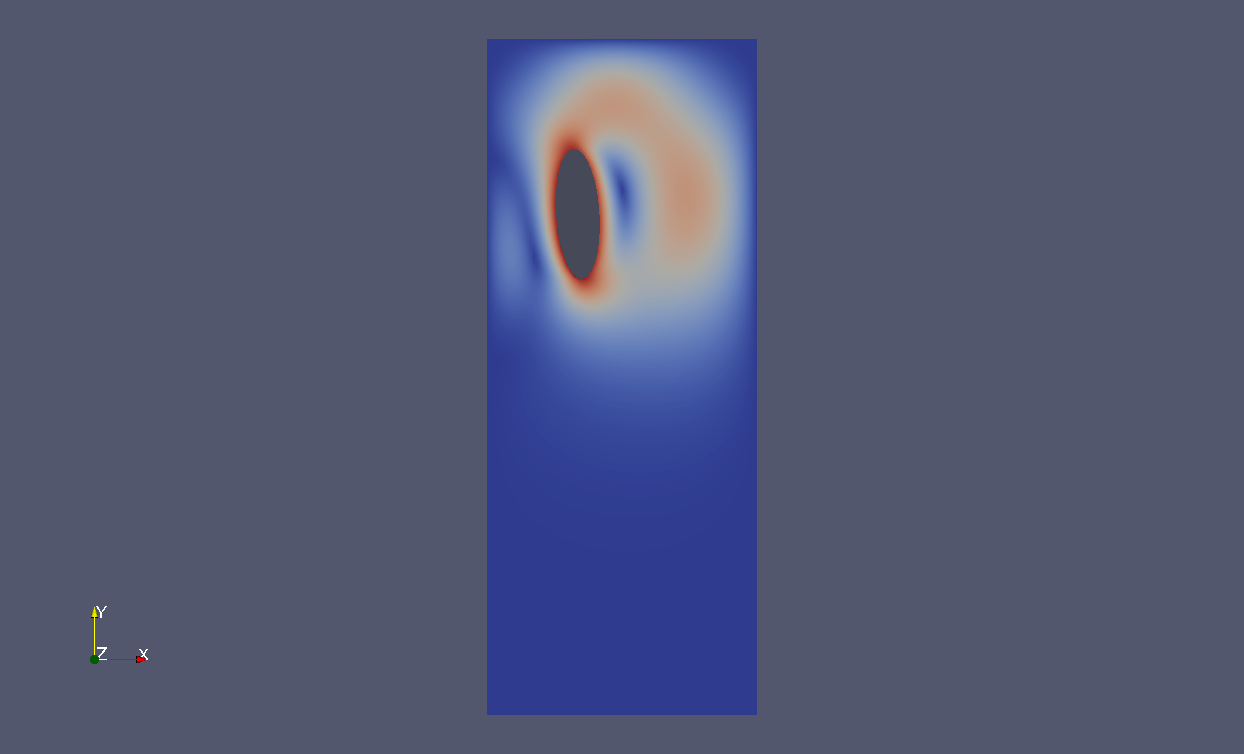}\\
 $t = 0$ & $t = 0.071$ & $t = 0.131$ & $t= 0.191$ & $ t= 0.251$\\
\end{tabular}
\end{center}

\begin{center}
\begin{tabular} {ccccc}
\includegraphics[trim = 16.5cm 1cm 16.5cm 1cm, clip, scale=0.25]{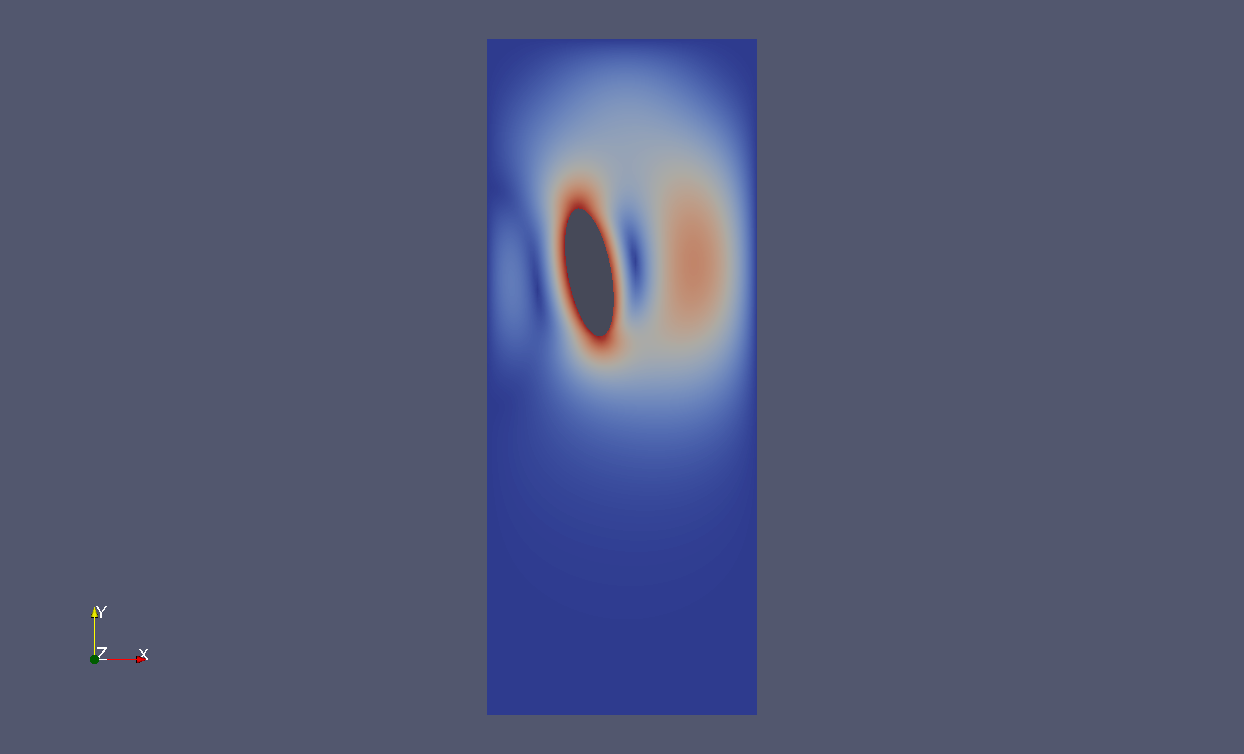} &
\includegraphics[trim = 16.5cm 1cm 16.5cm 1cm, clip, scale=0.25]{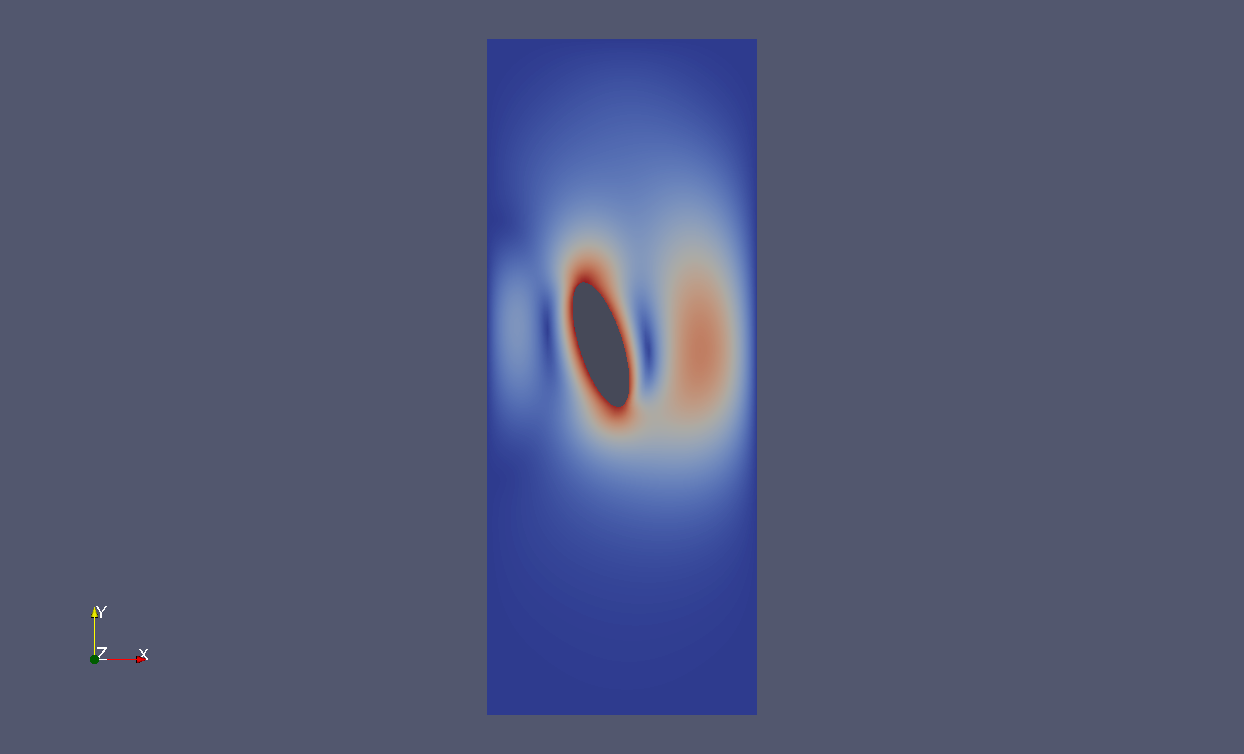} &
\includegraphics[trim = 16.5cm 1cm 16.5cm 1cm, clip, scale=0.25]{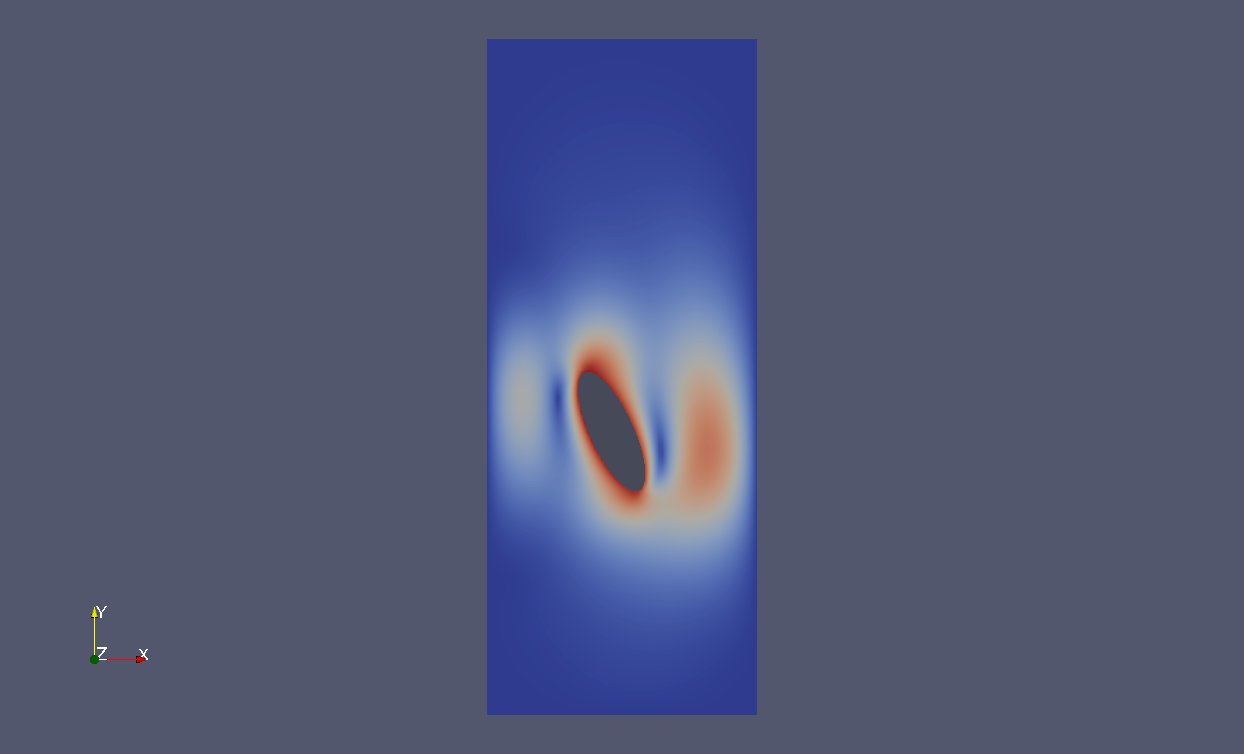} &
\includegraphics[trim = 16.5cm 1cm 16.5cm 1cm, clip, scale=0.25]{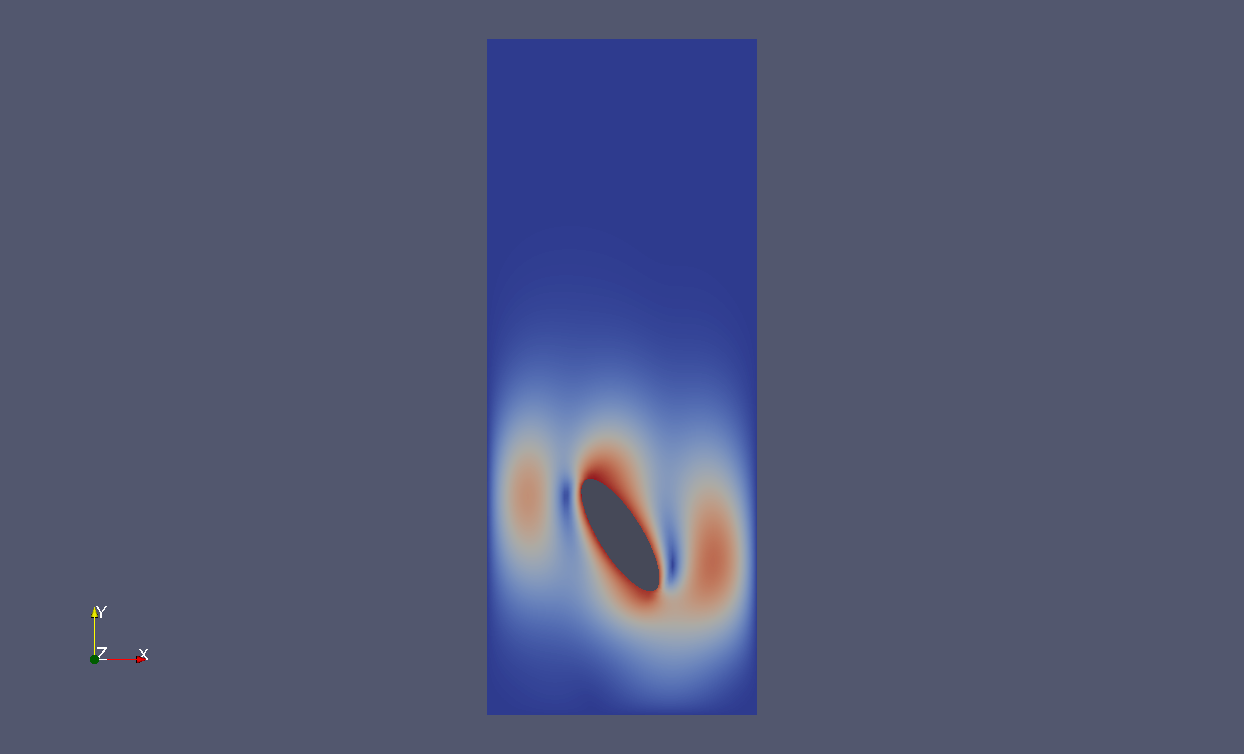} &
\includegraphics[trim = 16.5cm 1cm 16.5cm 1cm, clip, scale=0.25]{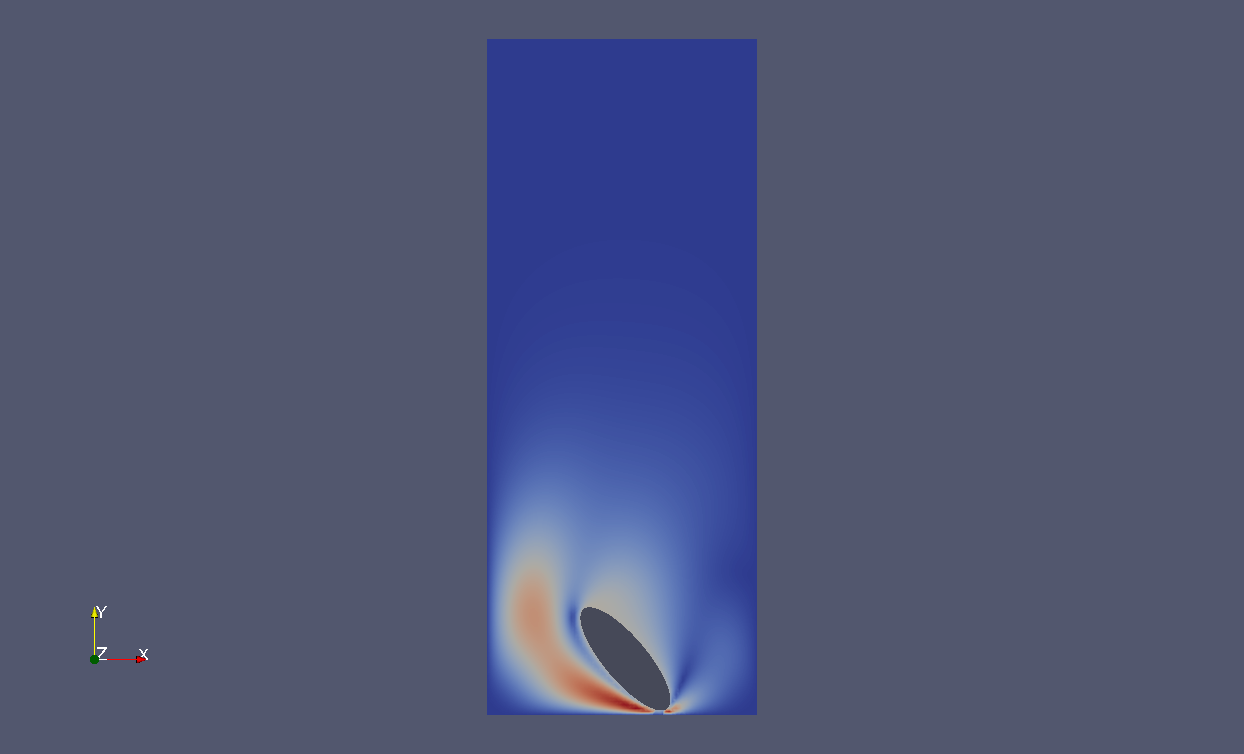}\\
$t = 311$ & $t = 371$ & $t = 431$  & $t = 491$  & $t = 551$\\
\end{tabular}\\
\begin{figure}[!h]
\centering
\caption{Fall of the solid.} \label{panorama}
\end{figure} 
\end{center}
In this simulation, we observe that the ellipse starts with straightening up and recentering in the channel, before turning over. Note that the good behavior of its dynamics would not be possible without in particular the stabilization technique which provides a good approximation of the forces exerted by the fluid.

\section{Conclusion}
In this work we have considered a fictitious domain method based on the ideas of Xfem, combined with a stabilization technique, and that we have applied to a Stokes problem and the Navier-Stokes equations coupled with a moving solid. The interest of this method lies mainly in two points: First the simplicity of the implementation, since all the variables (primal variables and multipliers) are defined on a single mesh which is independent of the computational domain. Secondly the robustness with regards to the computation of the normal trace of the Cauchy stress tensor whatever the way the computational domain cut the mesh. This second point is crucial in fluid-structure interaction models, because of the importance of the role played by the fluid forces. Applications to the simulation of the swim of deformable solids, applications in 3D and adaptation to other models constitute works in progress.

%\section*{Acknowledgments}
%This work is partially supported by the ANR-project CISIFS: 09-BLAN-0213-03, and the foundation STAE in the context of the RTRA platform SMARTWING. It is in particular based on the use of the Getfem++ library \cite{Getfem}, and thus on collaborative efforts with Yves Renard.

\bibliographystyle{plain}

\end{document}